\documentclass[12pt,a4paper]{article}
\usepackage[english]{babel}
\usepackage[latin1]{inputenc}
\usepackage[T1]{fontenc}
\usepackage{amsmath,amssymb}
\usepackage{graphicx}
\usepackage{caption}
\usepackage{subfig}
\usepackage{multirow}
\usepackage{array}
\usepackage{hyperref}
\usepackage{caption}
\usepackage{subfig}
\usepackage{multirow}
\usepackage{bm}
\usepackage{fancyhdr}
\usepackage{color}
\usepackage{changepage}
\usepackage{authblk}
{\left\lbrace\begin{array}{@{}l@{}}}%
{\end{array}\right.}
\title{Dynamical systems on graphs through the signless Laplacian matrix}
\author[1]{B. Giunti\thanks{barbara.giunti01@universitadipavia.it}}
\author[2]{V. Perri\thanks{vincenzo.perri538@edu.unito.it}}
\affil[1]{\footnotesize{Dipartimento di Matematica, Universit\`a degli Studi di Pavia, Pavia}}
\affil[1]{INFN, Sezione di Pavia, via Bassi 6, 27100 Pavia, Italy}
\affil[2]{\footnotesize{Dipartimento di Fisica dei sistemi complessi, Universit\`a degli Studi di Torino, Torino}}
\date{\vspace{-5ex}}
\begin{document}
\maketitle
\begin{abstract}
There is a deep and interesting connection between the topological properties of a graph and the behaviour of the dynamical system defined on it.
We analyse various kind of graphs, with different contrasting connectivity or degree characteristics, using the signless Laplacian matrix. 
We expose the theoretical results about the eigenvalue of the matrix and how they are related to the dynamical system. 
Then, we perform numerical computations on real-like graphs and observe the resulting system. 
Comparing the theoretical and numerical results we found a perfect consistency. 
Furthermore, we define a metric which takes in account the "rigidity" of the graph and enables us to relate all together the topological properties of the graph, the signless Laplacian matrix and the dynamical system.
\end{abstract}
\section{Introduction}
The interaction terms in many dynamical systems can be viewed as links in a graph, and it is interesting to explore new ways to study these systems by using graph theory. 
In our work, we discuss systems based on the signless Laplacian of a graph, motivated by models arising in the study of viral capsids \cite{CerInTwa}\cite{Giunti}\cite{Perri}.
The original idea of this model is due Cermelli, Indelicato and Zappa \cite{CerInZap}, who worked with a viral capsid, i.e., the "head" of a virus. 
Inside the capsid is contained the RNAm of the virus, and it opens in order to release its infective content. 
Therefore, those dynamical systems become relevant to study.
A rich literature on dynamical systems based on the Laplacian exists, for instance in the study of the synchronisation problem or for harmonic coupled oscillators. 
The problem may seem similar, but these two matrices have completely different behaviours. 
For more information about the Laplacian matrix and the synchronisation problem, see for example \cite{Syn01}\cite{Syn02}.  
The main difference between our case and the Laplacian is that the associated energy for the latter is not coercive: 
the asymptotic behaviour of the system involves large oscillations.
\\
Given a graph, we define a linear dynamical system on it. 
Each node represents a component of a vector field, while the links between nodes are the mutual influence in the potential function. 
To study the asymptotic behaviour of the system, we first characterise analytically the eigenvalues and the eigenvectors of the signless Laplacian for some general classes of graphs, given the explicit form of the eigenvalue and eigenvectors. 
Then, we perform numerical simulations on more complex systems, based on large and random graphs. 
In order to study the connection between the behaviour of the system and the topology of the graph we try, using specifically designed graphs, to isolate the various characteristics that seem meaningful, in order to understand the role of nodes with different properties. 
We found that not only the degree of a node is important, but also the average degree of its first, second and third neighbourhoods. 
We defined a metric to take this property into account. 
Using this distance, the simulations show a perfect match with the theoretical result and the expected behaviour.
The idea is that there is a kind of rigidity: the more the nodes are connected (the graph is "rigid"), the less they move. 
This property acts also locally: if just a few nodes are poorly connected, and the rest of the graph has a lot of links, they spread crazily while the other parts remain steady.
Then, we investigate also the importance of the clustering as a graph measure. 
It turns out to have some kind of significance. 
Therefore, we use the previous metric to define a new one, which is able to include the information of the clustering. 
Both measures are really valid to describe the behaviour of the system.
\section{The dynamical system}
Given a graph $G=(V,E)$ of order $n$ with vertex set $V$, edge set $E$ and adjacency matrix $A=(A_{ij})$, consider a vector field $\pmb{x}\in \mathbb{R}^{n}$ whose components $\left(x_{1},\dots,x_{n}\right)$ are scalar fields associated to the vertices of $G$.
\\
The signless Laplacian matrix of $G$ is $Q=A+D$, where 
$$D=\left(d_{ij}\right)=\delta_{ij} deg\left(v_{i}\right).$$
$Q$ is semidefinite positive, with non negative entries.
\\
We consider interactions among the nodes described by a potential function \cite{CerInTwa} $U: \mathbb{R}^{n}\rightarrow \mathbb{R}$ of the form 
\begin{align}
U\left(\pmb{x}\right)&=
\frac{a}{2}\sum_{i=1}^{n}x_{i}^{2}+\frac{b}{2}\sum_{i=1}^{n}\sum_{j=1}^{n}A_{ij}\left(x_{i}+x_{j}\right)^{2}
\\
&= 
\frac{a}{2}\sum_{i=1}^{n}x_{i}^{2}
+
\frac{b}{2}\sum_{i=1}^{n}\deg\left(i\right)x_{i}^{2}
+
\frac{b}{2}\sum_{j=1}^{n}\deg\left(j\right)x_{j}^{2}
+
b\sum_{i=1}^{n}\sum_{j=1}^{n}A_{ij}x_{i}x_{j}
\nonumber
\\
&=
\frac{1}{2}\sum_{i=1}^{n}\left(a+2b  \deg\left(i\right)\right)x_{i}^{2}+b\sum_{i=1}^{n}\sum_{j=1}^{n}A_{ij}x_{i}x_{j}
=\frac{1}{2}\pmb{x}\cdot H\pmb{x}
\nonumber
\end{align}
where $H=\left(a I + 2b Q\right)$, and $a,b$ are real parameters, with $b>0$ and, at least at the beginning, $a<0$. 
The associated gradient system has the form 
\begin{equation}
\dot{\pmb{x}}=-\nabla U \left(\pmb{x}\right)
\end{equation}
Note that the system is linear and $H$ symmetric, i.e., $H$ has got real eigenvalues. 
Denoting by $\lambda_{i}$ the eigenvalues of the dynamical system and $q_{i}$ the eigenvalues of the signless Laplacian $Q$, ordered in descending order, we have
\begin{equation}
\lambda_{i}=-a-2b q_{n-i+1}.
\label{relation}
\end{equation}
If the largest eigenvalue of the system is negative, the system is stable, so we have the simple stability condition:
$$
\text{if} \quad q_{n}>-\dfrac{a}{2b} \quad \text{and}\quad b>0
$$
then the system is stable, otherwise the system is unstable.
\section{$Q$-eigenvalues and $Q$-eigenvectors}
We compute the eigenvalues and the eigenvectors of the signless Laplacian for five classes of graphs: 
the complete graph on $n$ vertices $K_{n}$, the complete bipartite graphs $K_{n,m}$, the cycle $C_{n}$, the path $P_{n}$ and the star $S_{n}= K_{1, n-1}$.
The $Q$-eigenvalues will be always denoted by $q_{i}$, $i=1,\dots,n$, and  will be sorted in descending order $q_{1}\geq q_{2}\geq \dots \geq q_{n}$.
\\
$\omega_{k}$ will represent the $k$-th $n$-root of the unit. 
Recall that $\omega_{k}=\cos\left(\frac{2\pi k}{n}\right)+i \sin\left(\frac{2\pi k}{n}\right)$. 
Further, we define 
$U_{n}=\left\lbrace \pmb{v} \in \mathbb{R}^{n}\mid \sum_{i=1}^{n}v_{i}=0\right\rbrace$.
The table reports the $Q$-eigenvalues and the $Q$-eigenvectors for each of the five classes of graphs.
\begin{figure}[h!]
\begin{adjustwidth}{-2,7cm}{}
\begin{tabular}{|| c | c | c | c || }
\hline\hline
Graph & $Q$-eigenvalues & Multiplicity & $Q$-eigenvectors
\\
\hline\hline
$K_{n}$ & $q_{1}=2n-2$ & $1$ & $\pmb{v}_{1}=\left(1,\dots,1\right)$ 
\\  
& $q_{n}= n-2$ & $n-1$ & $\sum_{i=1}^{n}v_{n_{i}}=0$
\\ 
\hline
$K_{n,m}$ & $q_{1}=n+m$ & $1$ & $\pmb{v}_{1}=\left(\frac{1}{n},\dots,\frac{1}{n}, \frac{1}{m}\dots,\frac{1}{m}\right)$ 
\\  
& $q_{2}= n$ & $m-1$ & $ \pmb{v}_{2}\in U_{n} \oplus \left\lbrace 0\right\rbrace^{m}$
\\
& $q_{m}= m$ & $n-1$ & $ \pmb{v}_{m}\in \left\lbrace 0\right\rbrace^{n}\oplus U_{m}$
\\
& $q_{n+m}= 0$ & $1$ & $ \pmb{v}_{n+m}=\left(1,\dots,1, -1\dots,-1\right)$ 
\\
\hline
$C_{n}$ & & &
\\
$n$ even & $q_{1}=4$ & $1$ & $\pmb{v}_{1}=\left(1,\dots,1\right)$ 
\\  
$i=1, \dots , \dfrac{n}{2}-1$ & $q_{2i}= 2+2\cos\left(\dfrac{2\pi i}{n}\right)$ & $2$ & $\pmb{v}_{2i}=
\left(
1, \cos\left(\dfrac{2\pi i}{n}\right),
\dots, \cos\left(\left(n-1\right)\dfrac{2\pi i}{n}\right)
\right)$
\\
& & & $\pmb{v}_{2i}^{'}=
\left(
0, \sin\left(\dfrac{2\pi i}{n}\right),
\dots, \sin\left(\left(n-1\right)\dfrac{2\pi i}{n}\right)
\right)$
\\
& $q_{n}=0$ & $1$ & $\pmb{v}_{1}=\left(1,-1,\dots,1,-1\right)$ 
\\
& & &
\\
$n$ odd & $q_{1}=4$ & $1$ & $\pmb{v}_{1}=\left(1,\dots,1\right)$ 
\\  
$i=1, \dots , \dfrac{n-1}{2}$ & $q_{2i}= 2+2\cos\left(\dfrac{2\pi i}{n}\right)$ & $2$ & $\pmb{v}_{2i}=
\left(
1, \cos\left(\dfrac{2\pi i}{n}\right),
\dots, \cos\left(\left(n-1\right)\dfrac{2\pi i}{n}\right)
\right)$
\\
& & & $\pmb{v}_{2i}^{'}=
\left(
0, \sin\left(\dfrac{2\pi i}{n}\right),
\dots, \sin\left(\left(n-1\right)\dfrac{2\pi i}{n}\right)
\right)$
\\
\hline
$P_{n}$ & $q_{n}=0$ & $1$ & $\pmb{v}_{1}=\left(1,\dots,1\right)$ 
\\
$i=1,\dots, n-1$ & $q_{i}= 2-2\cos\left(\dfrac{2\pi i}{n}\right)$ & $1$ & 
$\pmb{v}_{i}=
\left(2+2\cos\left(\dfrac{2\pi i}{n}\right),\dots, \right.$
\\
& & & $2\cos\left(k\dfrac{2\pi i}{n}\right)+2\cos\left(k+1\dfrac{2\pi i}{n}\right), $
\\
& & & $\left.\dots, 2+2\cos\left(\dfrac{2\pi i}{n}\right)\right)$
\\
\hline
$S_{n}$ & $q_{1}=n$ & $1$ & $\pmb{v}_{1}=\left(n-1, 1, \dots, 1 \right)$
\\
& $q_{2}=1$ & $n-2$ & $\displaystyle \sum_{i=1}^{n}v_{i}=0$
\\
& $q_{n}=0$ & $1$ & $\pmb{v}_{1}=\left(-1, 1, \dots, 1 \right)$
\\
\hline\hline
\end{tabular}
\end{adjustwidth}
\end{figure}
We can see how the multiplicity of the $Q$-eigenvalue decreases as the graphs become less and less connected.
\\ 
Moreover, there is a substantial difference for $C_{n}$ if $n$ is odd or even. 
The reason lies in the bipartiteness of $C_{n}$ for $n$ even. 
In fact, $C_{2k}$ is bipartite while $C_{2k+1}$ is not. 
As a result, the smallest $Q$-eigenvalue of $C_{2k}$ is unique and equal to $0$\cite{CveSiI}. 
On the other hand, $C_{2k+1}$ as a smallest $Q$-eigenvalue of multiplicity $2$ and greater than $0$. 
As a consequence, the dynamical system based on $C_{2k}$ will turn out to be more likely unstable.
In order to avoid this problem, we will no longer impose the restriction $a<0$.
\section{Dynamical system}
Due to the relation (\ref{relation}) between the eigenvalues of the dynamical system and the $Q$-eigenvalues, we want some upper and bound for the largest and the smallest $Q$-eigenvalues. 
There are several results in regard to, for a general review see \cite{CveSiI} or \cite{Sta}.
We discuss here only the two systems based on $K_{n}$ and on $P_{n}$, in order to point out the general trend. 
The $K_{n}$ system has smallest eigenvalue $\lambda_{n}=-a-4b\left(n-1\right)$, with eigenvector $\left(1,\dots,1\right)$, which corresponds to a so-called 'breathing' mode, i.e., a perturbation where all vertices spread together at the same speed.  
As to the other instability modes, i.e., the other eigenvectors of the dynamical system, we have that all other $\left(n-1\right)$ $Q$-eigenvalues of $K_{n}$ are equal to $\left(n-2\right)$ and so all other eigenvalues of $K_{n}$ are $\lambda_{i}=-a-2b\left(n-2\right)$ for $i=1,\dots, n-1$.
The corresponding eigenvectors $\pmb{v}_{i}$, for $i=1,\dots, n-1$, are such that 
$$\sum_{i=1}^{n}\pmb{v}_{i}=0.$$
We analyse the other eigenvalues: 
if $a<0$, then $-a$ is positive while $-2b\left(n-2\right)$ is negative, more and more negative growing $n$. 
That means it should be easy to make the system stable because growing $n$ to destabilise the system we need to increase $a$ very fast. 
Though, for $n>>1$ if $b>1$ also for great value of $a$ the system will be stable.
If $a>0$, then $-a$ is negative and so it is $-2b\left(n-2\right)$, and the system is in general stable.
\\
\\
All eigenvalues of $P_{n}$ have multiplicity equal to one. 
Thus, all the eigenspaces have the dimension equal to one and each perturbation goes along just one vector. 
\\
The path is a bipartite graph and so the smallest $Q$-eigenvalue is equal to zero \cite{CveSiI}. 
Then, the system's largest eigenvalue is $\lambda_{1}=-a$, with eigenvector $\pmb{v}\left(q_{n}\right)=\left(1,\dots,1\right)$.
When we perturb the system along this eigenvector we find again a breathing perturbation:
$$ 
\pmb{x}\left(t\right)=e^{\lambda t}c_{1}\left(1,\dots,1\right)=\left(c_{1}e^{\lambda t},\dots,c_{1}e^{\lambda t}\right). 
$$
The other system's eigenvalues are $\lambda_{j}=-a-4b\left(1-\cos\left(\frac{2 \pi j}{n}\right)\right)$, with corresponding eigenvectors:
\begin{equation}
\begin{split}
\pmb{v}\left(\omega_{j}\right)
=&
\left(
 2+2\cos\left(\frac{2 \pi j}{n}\right), \dots, 2\cos\left(\frac{2 \pi j k}{n}\right)+2\cos\left(\frac{2 \pi j \left( k+1 \right)}{n}\right),\right. \\
&\quad
\left.
\dots,
2+2\cos\left(\frac{2 \pi j}{n}\right)
\right).
\end{split}
\end{equation}
It is interesting to notice that this system is really easy to perturb, i.e., a weak force in one direction can disrupt it.
Although, the perturbation acts in that direction only, leaving the other vertices untouched. 
\section{Simulations on large graphs}
We will now analyse some numerical results. 
We want to move on from the previous theoretical results and to test if there is some property of the graph that directly influence the stability of the system. 
Then, we study the asymptotic convergence of the system just beyond the stability threshold.
\\
Recall that the eigenvalues of $Q$ are all non-negative and related to the system's ones by the relation (\ref{relation}). 
The stability diagram of the system is given plotting for each eigenvalue the values of $a$ and $b$ for which $\lambda_i=0$.
\begin{figure}[h]
\centering 
\begin{tabular}{cc}
\includegraphics[width=.5\columnwidth]{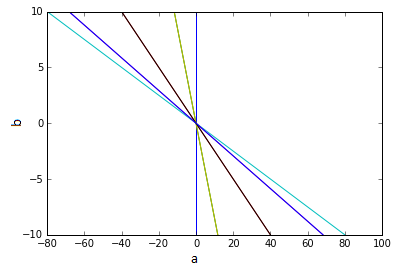} 
\includegraphics[width=.5\columnwidth]{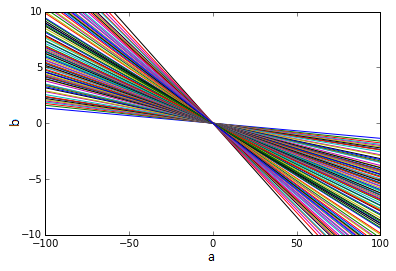}
\end{tabular}
\caption{$(a)$ Stability diagram of a cycle graph with n=8. $(b)$ Stability diagram for a scale-free system with n=100, obtained using the Barab\'asi and Albert model. Differently from $(a)$ in $(b)$ the first area beyond the stability threshold is not easily identifiable. }
\label{stab}
\end{figure}
We are interested in studying what happens when $a$ and $b$ have values between the straight lines of the first and the second eigenvalue. In other words, when the largest eigenvalue is positive while all the others remain negative (see Figure \ref{stab}).
Using (\ref{relation}) for the largest eigenvalues, we impose the values of $a \in \mathbb{R}$ and $b \in \mathbb{R}^+$ such that  $\lambda_1>0$ and $\lambda_2<0$:
\begin{equation}
\left.
\begin{array}{r}
\lambda_{1}= -a-2bq_{n}>0 \ \ \rightarrow \ \ \  a<-2bq_{n}  
\\
\lambda_{2}= -a-2bq_{n-1}<0 \rightarrow  a>-2bq_{n-1} 
\end{array}
\right\rbrace
\Rightarrow
-2bq_{n-1}<a<-2bq_{n}
\end{equation}
The general solution of the system contains as exponents the eigenvalues of the system. 
Then, imposing these conditions, all the contribution of the various eigenspaces should go asymptotically to zero, except the component of the largest eigenvalue.
Calling $E_i$ the eigenspace of the $i-$th eigenvalue and $\bm{x}(t)$ the state of the system at time $t$, we obtain 
\begin{equation}
dist\left(\bm{x}\left(t\right), E_{1}\right)\to 0 \quad \text{for}\quad t \to + \infty
\label{conv}
\end{equation}
if we do not have an initial condition $\bm{x}_0 \in E_i$, with $i \neq 1$. 
If as an initial condition we have a vector $\bm{x}_0 \in E_i$, $i \neq 1$, the system will asymptotically go to zero along the direction of the eigenspace. 
In the following, for convenience, we will often refer to the eigenspace $E_1$ and to its eigenvector and eigenvalue using the adjective \emph{principal}.
\begin{figure}[!ht]
\centering
\begin{tabular}{cc}
\subfloat[]{\includegraphics[width=.5\columnwidth]{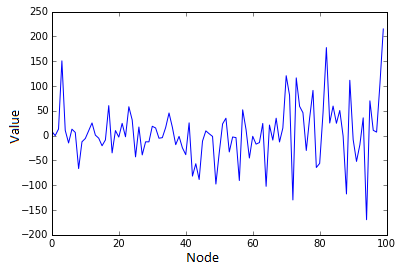}}
\subfloat[]{\includegraphics[width=.5\columnwidth]{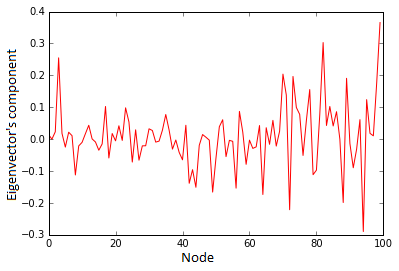}} 
\end{tabular}
\caption{$(a)$ The state of the system at an advanced numerical step. 
$(b)$ The dominant eigenvector of the system.}
\label{BA}
\end{figure}  
The computation shows the expected asymptotic convergence of the state of the system to $E_1$ (see Figure \ref{BA}).
\\
Given that the system converges to a state proportional to $E_1$ we want to study how
the shape of this state is related to the structural features of the graph
\\
Then, we perform a graphical comparison between the value of the component of the node in $E_i$ and different graph measures for the centrality of a node.
\begin{figure}[!ht]
\centering
\begin{tabular}{cc}
\subfloat[]{\includegraphics[width=.5\columnwidth]{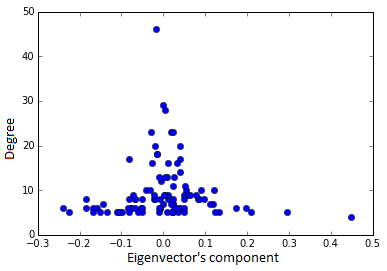}}
\subfloat[]{\includegraphics[width=.5\columnwidth]{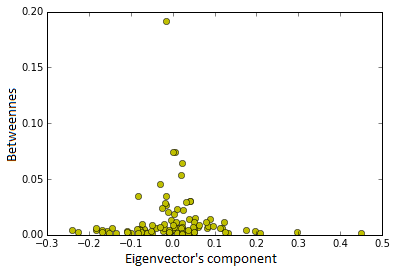}} 
\\
\subfloat[]{\includegraphics[width=.5\columnwidth]{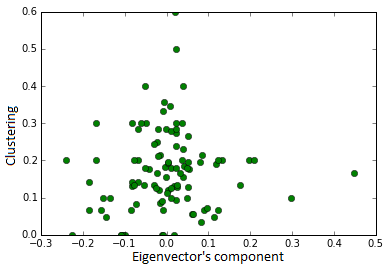}} 
\subfloat[]{\includegraphics[width=.5\columnwidth]{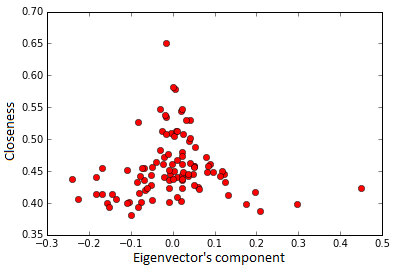}}  
\end{tabular}
\caption{ 
Comparison between the value of the node in the eigenvector component and its centrality in a graph generated using the Barab\'asi and Albert model. 
$(a)$ Degree. 
$(b)$ Betweenness.  
$(c)$ Clustering. 
$(d)$ Closeness. 
In all the centrality measures but the clustering a tendency for isolated nodes to have a higher eigenvector component is observed.}
\label{BA1}
\end{figure}
We can see a marked tendency for isolated nodes to have a higher weight in $E_i$. 
However, from this observation, we cannot affirm that the perturbation accumulates on the isolated nodes because this is not true for every node with a low centrality. 
For example (Figure \ref{BA1}), isolated nodes with a low value in the eigenvector, i.e., nodes on which the perturbation does not act, can be easily found. 
Therefore, we can only conclude that if a node has a large component in the eigenvector then it is an isolated one and that highly central nodes have low values in the eigenvector.
\section{Detailed study of isolated nodes}
\begin{figure}[!ht]
\centering
\includegraphics[width=.6\columnwidth]{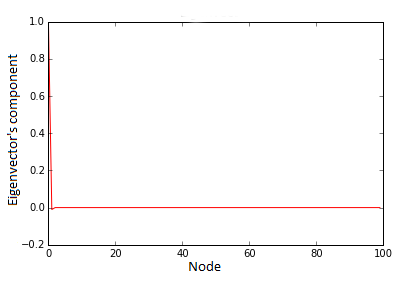}
\caption{Eigenvector for the graph obtained connecting a node to a full graph, with one edge only.}
\label{ci}
\end{figure}

We need to analyse which characteristic of an isolated node give it a high value in the eigenvector. 
Firstly we want to test the behaviour of a graph that is highly connected except for a single dangling node, (i.e., a node with degree one). 
We connect an isolated node to a $n=100$ complete graph. 
The isolate node's value in the eigenvector is $0.99$ which is almost the highest value it could reach since the eigenvector's module is normalised to one (see Figure \ref{ci}).
Only another node has a value in the eigenvector, $-1.0\cdot10^{-2}$, that is significantly greater than zero, while all the other nodes have an extremely small value: $5.3\cdot10^{-5}$.
That is the node connected to the isolated one.
This is also the node with the highest degree.
That suggests that the neighbours of a node have a role determining its value inside the eigenvector. 
Therefore not only the
centrality of a node is important but also the ones of its neighbours 
\\
To investigate how it works, we generate a scale-free graph with $n=100$ and modify it adding one or two nodes in three ways: 
one to the node with the highest degree, or
one to the node with the lowest degree, or
one to the node with the highest degree and the other to the lowest degree one.
\begin{figure}[!ht]
\centering
\begin{tabular}{ccc}
\subfloat[]{\includegraphics[width=.33\columnwidth]{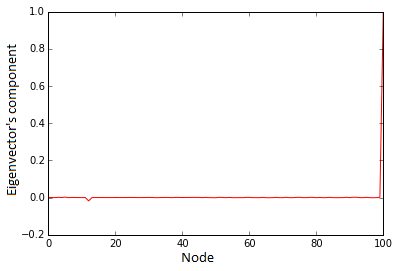}}
\subfloat[]{\includegraphics[width=.33\columnwidth]{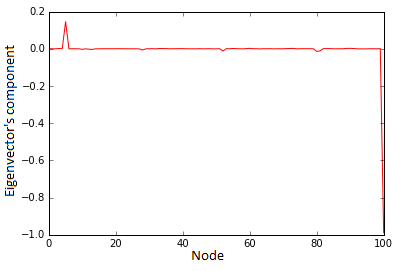}} 
\subfloat[]{\includegraphics[width=.33\columnwidth]{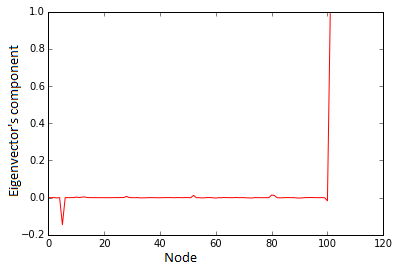}} 
\end{tabular}
\caption{Comparison between the eigenvectors we obtain adding one or two nodes to a scale-free graph. $(a)$ One node added and connected to the node with the highest degree. $(b)$ One node added and connected to the node with the lowest degree.  $(c)$ Two nodes added and connected one to the node with the highest degree and the other to the lowest degree one}
\label{collegato}
\end{figure} 
We analyse the resulting eigenvectors.
\\
In the first case, the dangling node connected to the hub ($(a)$ in Figure \ref{collegato}) takes the maximum value in the eigenvector: $9.9\cdot10^{-1}$.
The hub is the only other node taking a value distinguishable from zero: $-1.8\cdot10^{-2}$. 
In the second case, we still see the isolated node taking the highest value $-9.8\cdot10^{-1}$.
Moreover, the node to which it is connected has the value $1.45\cdot10^{-1}$, that is almost ten times larger than the one of the hubs in the previous case.
Furthermore, also other nodes have a value which is distinguishable from zero.
These are the lower degree neighbours of the pendant node.
In the third case ($(c)$ in Figure \ref{collegato}), the eigenvector is almost identical to the second case.
The nodes show almost the same values in the eigenvectors (even if with reversed sign), with the only exception of the dangling node connected to the hub  which has a rather low value, $-1.7\cdot10^{-2}$.
Note that this value is on the order of the values shown by the lower degree neighbours of the lowest degree node. It means only one of the added nodes take a remarkable value, while the other barely moves.
\\
These results underline the importance of neighbours and show the presence of some kind of "winner takes all" mechanism. 
In fact, in the third graph where two nodes are added, only one of the dangling nodes takes all the value.
\section{The importance of neighbours}
\begin{figure}[!ht]
\centering
\includegraphics[width=.5\columnwidth]{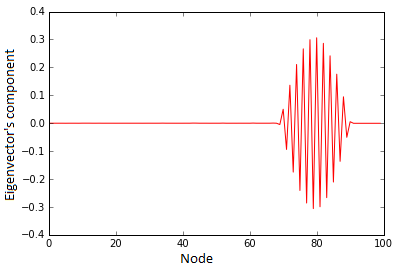}
\caption{Principal eigenvector of a path that connects a Barab\'asi Albert and a complete graph.  }
\label{ultimi}
\end{figure}

\begin{figure}[!ht]
\centering
\begin{tabular}{cc}
\subfloat[]{\includegraphics[width=.5\columnwidth]{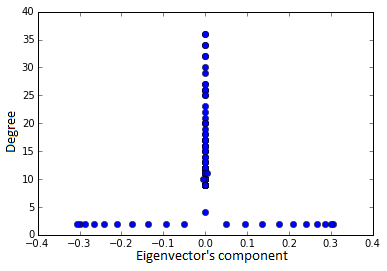}}
\subfloat[]{\includegraphics[width=.5\columnwidth]{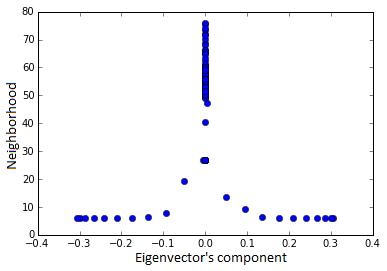}}
\end{tabular}
\caption{Comparison between centrality measures and the principal eigenvector for the graph generated connecting Barab\'asi Albert to a complete graph through a path .$(a)$ Degree. $(b)$ Closeness. $(d)$ The $x$ introduced in \eqref{c1}.}
\label{ultimi2}
\end{figure}

\begin{figure}[!ht]
\centering
\begin{tabular}{cc}
\subfloat[]{\includegraphics[width=.5\columnwidth]{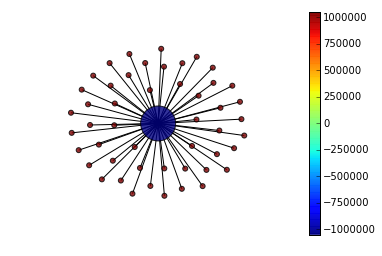}}
\subfloat[]{\includegraphics[width=.5\columnwidth]{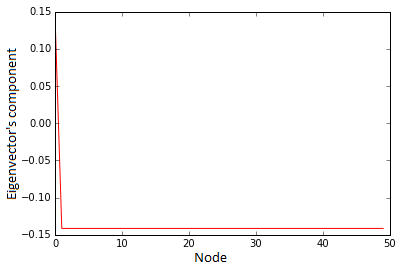}}\\
\end{tabular}
\caption{$(a)$ Star graph. $(b)$ Principal eigenvector of the star graph.}
\label{Starg}
\end{figure}

\begin{figure}[!htb]
\centering
\begin{tabular}{cc}
\subfloat[]{\includegraphics[width=.5\columnwidth]{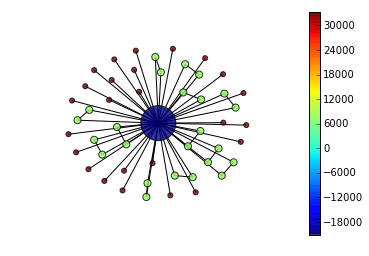}}
\subfloat[]{\includegraphics[width=.5\columnwidth]{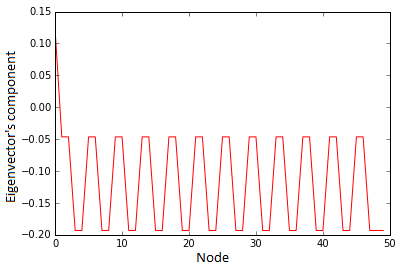}}\\
\end{tabular}
\caption{$(a)$ Modified Star graph. $(b)$ Principal eigenvector of the Modified star graph.}
\label{Star}
\end{figure}

We noted that a notion of rigidity could be involved as a part of the explanation for the observed behaviour. 
Recall the previous result on the graph, especially the comparison between the $K_{n}$ system and the $P_{n}$ system. 
In the $K_n$ system, we see a rise in stability with the growth of $n$. 
This can be related to the number of neighbours that a node has got: a higher number implies a higher constraint. 
Also, the numerical results show the same trend.
For example, in Figure \ref{collegato}, the connection of a dangling node to a hub limits the propagation of a perturbation. 
Moreover, the hub itself has a low value in the eigenvector.
On the other hand, the connection to a less constrained node allows the perturbation to propagate to its lower degree neighbours.
To analyse this hypothesis, we create a graph in which areas with different rigidity are easily recognisable. 
This graph is built connecting a path of twenty nodes to a scale-free graph of seventy nodes and to a complete graph of ten nodes. 
In these graphs, it is clear that the path is the less rigid area.
Thus, our hypothesis suggests that the non-zero values of the eigenvector should accumulate there. 
\\
Computing the eigenvectors we find that both rigid parts have null components in the eigenvector. 
Instead, the components of the nodes in the path are different from zero. 
The sign for each node is opposite than the ones of its neighbours and their magnitude increase moving away from the rigid parts. 
The maximum value is then reached by the nodes at the centre of the path (see Figure \ref{ultimi}).   
This result confirms the assumption since the larger components are those of the path. 
It can also be noticed that both the complete and the scale-free graphs have null components.\\
Taking this into account, we now want to build a measure that can capture the rigidity of the nodes. 
We say that the rigidity of a node comes from its degree and from the constriction imposed by the presence of its neighbours.
Hence, we consider the medium degree of the first, second and third neighbours of a node and add their value (weighted, in order to consider the lower influence of a distant node) to the node's degree 
\begin{equation}
\begin{split}
r\left( i \right)
:=
&
\deg(i)+\frac{p}{|N(i)|}\sum_{j\in N(i)} \deg(j)
+
\\
&
+
\frac{p^2}{| N(N(i))|}\sum_{k\in N(N(i))} \deg(k)
+
\frac{p^3}{|N(N(N(i)))|}\sum_{u\in N(N(N(i)))} \deg(u)
\end{split}
\label{c1}
\end{equation}
where $0\le p\le 1$ indicates the influence that the neighbours' degree of the node $i$ has over its rigidity and $N(i)$ is the set of neighbours of node $i$. 
In this formulation, the influence of a node at distance $k=1,2,3$ is reduced by a factor $p^k$.
Comparing the principal eigenvector with the centrality measure, we see in Figure \ref{ultimi2} that the degree alone does not distinguish between the various nodes in the path since they have all the same degree.
Considering also the neighbours' degree, the nodes of the path are subdivided giving to the ones with a lower value in the eigenvector a higher centrality.
Comparing the principal eigenvector with the centrality measures for the second graph we obtain the results in Figure \ref{ultimi2}.
This new rigidity measures seem to be able to capture some of the aspects that control the behaviour of the system, i.e., the importance of isolated nodes and the importance of the neighbours.
\\
From the results obtained above, we do not expect the form of the eigenvector to depend on the clustering of the nodes. In fact, in Figure \ref{BA1}, we see no relation between the importance of a node in the eigenvector and its clustering.
Nevertheless, analysing the original potential
\begin{equation*}
U\left(\bm{x}\right)=\frac{a}{2}\sum_{i=1}^{n}x_{i}^{2}+\frac{b}{2}\sum_{i=1}^{n}\sum_{j=1}^{n}A_{ij}\left(x_{i}+x_{j}\right)^{2}
\end{equation*}
we consider the interaction term for two nodes only and search in which direction it is minimised. 
Excluding the solution $\bm{x}=\bm{0}$ since the parameters $a$ and $b$ are chosen to make it unstable, we have that the minimum is given by configuration $x_i=-x_j$. 
This explains why neighbouring nodes tend to have values with opposite sign (particularly evident in Figure \ref{ultimi2}). 
This means that the neighbouring nodes of the system assume values with opposite sign trying to minimise the energy. 
\\
A perfect example is the star graph (see Figure \ref{Starg}).
In this graph, the nodes can easily minimise their energy assuming a value whose sign is the opposite of the one of the neighbour because each node but the central one has got just one neighbour: the central node itself.
\\
However, if in $S_{50}$ we introduce an edge between some of the nodes, it is impossible to minimise the energy assuming opposite values. 
The system tries to overcome this obstacle giving these nodes a low value (see Figure \ref{Star}). 
The node at the centre of the star, however, cannot assume a low value since it has to minimise the energy of the interaction also with the other non-constrained nodes. 
Therefore, forced between these two requests, it assumes an intermediate value.
\begin{figure}[!htb]
\centering
\begin{tabular}{cc}
\subfloat[]{\includegraphics[width=.45\columnwidth]{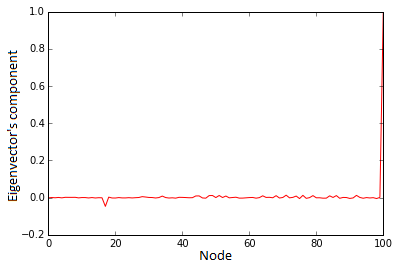}}&\subfloat[]{\includegraphics[width=.45\columnwidth]{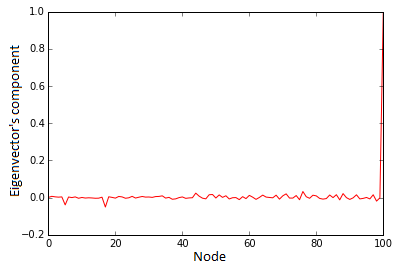}}
\\
\subfloat[]{\includegraphics[width=.45\columnwidth]{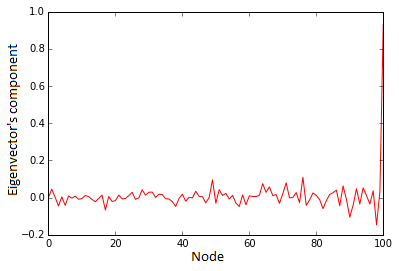}}&\subfloat[]{\includegraphics[width=.45\columnwidth]{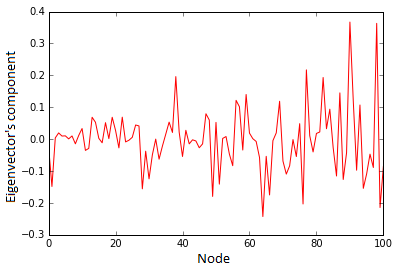}}
\\
\end{tabular}
\caption{Barab\'asi Albert graph having a principal node with high clustering. A node is added and connected to $(a)$  to the highest degree hub $(b)$ the two highest degree hubs  $(c)$ the three highest degree hubs and $(d)$ the four highest degree hubs.}
\label{NC}
\end{figure} 
Since the clustering of a node is the measure of the connections between the neighbours of a node, we expect it to have some role determining the value of the component in the principal eigenvector. 
On the contrary, it seemed to have no connection with the principal eigenvector (see Figure \ref{BA1}). 
Then, for most of our analysis, this measure was calculated but rarely taken into account.
To investigate the role of the clustering, we use the BA algorithm to generate a new graph with $n=100$. The lowest value for the degree in this graph is five.  
Then we add a new node to the graph and generate four new graphs connecting the newly introduced node to the four higher degree nodes of the network.
The four hubs in the initial graph are connected together.
Thus, all the neighbours of the $101^{st}$ node are neighbours to each other, and so its clustering is $c(101)=1$. 
Note that, when it is connected to only one node, the clustering value is zero by default.
In Figures \ref{NC} $(a)$, $(b)$ and $(c)$, we see that the node has the highest value in the eigenvector. 
On the other hand, in case $(d)$ the value in the eigenvector is almost negligible $\left(9\times 10^{-2}\right)$, even though it still has the lowest degree.
Then, also the clustering plays a certain role.
However, for large graphs without trivial connection patterns, its role is marginal compared to that of the degree.
From this, we can suppose that the clustering can be used to improve the rigidity measure we defined earlier. Then, we define a new measure $\tilde{r}$:
\begin{equation}
\tilde{r}\left(i\right)
:=
r\left(i\right)
+
Clust\left(i\right)\tilde{p}\deg\left(i\right)
\end{equation}
where $Clust\left(i\right)$ and $\deg\left(i\right)$ are respectively the clustering coefficient and the degree of node $i$, while $\tilde{p}$ is an adjustable parameter. 
In this formula $\deg\left(i\right)$ multiplies $Clust\left(i\right)$ to differentiate between nodes that despite belonging to equally clustered groups does not have the same rigidity due to the different amount of links.
\begin{figure}[!h]
\centering
\begin{tabular}{cc}
\subfloat[]{\includegraphics[width=.5\columnwidth]{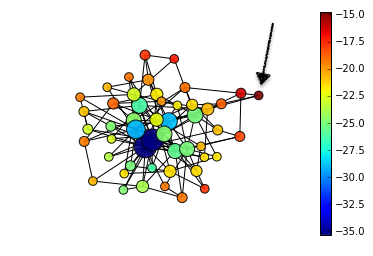}}
\subfloat[]{\includegraphics[width=.5\columnwidth]{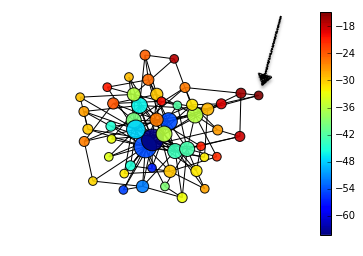}}\\
\end{tabular}
\caption{$(a)$ Random graph with the metric $r$. $(b)$ Random graph with the metric $\tilde{r}$.}
\label{rrtilde}
\end{figure} 
Using this new metric, we see (Figure \ref{rrtilde}) an improvement in the behaviour prediction. 
As before, the colour shows the value of the measure used, and the size of a node represents the degree of the node.
The arrow shows the node whose value was higher in the eigenvector. 
We can see that both the metric are able to capture this fact. Moreover, there is a tiny improvement using $\tilde{r}$ instead of $r$. This is due to the information carried by the clustering, as we said above.
\\
Summing up, those metric are both truly useful for studying the properties of the graph and of the system because they are able to discern which part of a graph are more able to move or spread and which not.  
\\
\section{Conclusions}
In the first part of the work, we analysed the stability of a dynamical system defined on several classes of graphs. The general result is that the more connected the graph is, the stronger under perturbation the system is. Keeping this in mind, in the second part of the work we tried to understand how exactly the connection properties influence the system, and also how we can measure of this connection properties. 
\\
From the theoretical results, we know that a perturbation spreads in different ways on the path and on the complete graph, i.e., on graphs with general low connectivity and on highly connected graphs. 
More precisely, if a node is connected with several others, and also those have a high connection rank, the perturbation is not going to have a huge impact. 
On the other hand, if the node is not very well connected with other nodes, or if it is connected to a node which degree is low, then the perturbation is going to have great effects. When we performed the numerical simulations, we found a perfect accord with this expectation. 
In fact, through numerical simulations, we proved that if we place the system in a specific type of instability its state converges in the direction of the principal eigenvector. 
This encouraged us to study the connection between  the features of this state and the structure of the graph. 
Performing simulations on different types of large graphs, constructed in order to have specific characteristics, we showed that exists a connection between the topology of the nodes and their importance inside the principal eigenvector. 
Our results give high relevance to the degree of a node, accompanied by the importance of the degree of its neighbours and the connections between them, carried out by the clustering.
\\
In this study is given a good account of the factors that lead a node to have a high value in the eigenvector. 
The relationship between the topology of the graph and its value inside the eigenvector, however, remains unknown. 
An analytical result that allows us to understand  the importance of a node from the graph structure knowing its topological characteristics, it was not obtained. 
This probably depends on the fact that the behaviour of a complex graph can not be simplified over a certain limit, as we tried to do connecting linearly the topological properties of the graph with the spectral ones.
\\
Since the behaviour of the system derives from a complex interplay of different characteristics, there is probably a more complex type of connection that has to be taken into account when we analyse the graph.
There is probably a more complex type of connection that has to be taken into account when we analyse the graph. 
The behaviour of the system derives from a complex interplay of the different characteristics of the graph.
\\
\\
\section*{Acknowledgement}
We would like to express our gratitude to Prof. Paolo Cermelli, who supported and encouraged us along all this work, from when it started as two Master thesis until this last version as a paper. 
Really thanks for being so helpful and still let us work on our on. It taught us a lot, and there is no greater compliment for a professor.

\end{document}